\documentclass{article}
\usepackage{amssymb}
\usepackage{amsfonts}
\usepackage{amsmath}
\usepackage{amsthm}
\usepackage[all]{xy}
\usepackage{ifthen}
\usepackage{makeidx}
\usepackage{enumerate}
\usepackage{eucal}
\usepackage{listings}
\usepackage{color}
\usepackage{fix-cm}
\usepackage{graphicx}
\usepackage{setspace}
\theoremstyle{theorem}
\newtheorem{theeorem}{Theorem~\!\!}
\newtheorem{theorem}{Theorem}

\newtheorem{example}{Example}

\begin{document}
\title{Effective Results on non-Archimedean Tropical Discriminants}

\author{Korben Rusek
\thanks{Partially supported by NSF MCS grant DMS-0915245.  
K.R.\ also partially supported by Sandia National
Labs and DOE ASCR grant DE-SC0002505. Sandia is a multiprogram
laboratory operated by Sandia Corp., a Lockheed Martin Company, for
the US DOE under Contract DE-AC04-94AL85000.} 
}

\date{\today}
\newcommand{\lcm}{\textrm{lcm}}
\newcommand{\irr}{\textrm{irr}}
\newcommand{\sylp}{$Syl_{p}$}
\newcommand{\BA}{\mathbb{A}}
\newcommand{\BN}{\mathbb{N}}
\newcommand{\BZ}{\mathbb{Z}}
\newcommand{\BQ}{\mathbb{Q}}
\newcommand{\BR}{\mathbb{R}}
\newcommand{\BC}{\mathbb{C}}
\newcommand{\BF}{\mathbb{F}}
\newcommand{\BP}{\mathbb{P}}
\newcommand{\BV}{\mathbb{V}}
\newcommand{\CA}{\mathcal{A}}
\newcommand{\CF}{\mathcal{F}}
\newcommand{\CQ}{\mathcal{Q}}
\newcommand{\fa}{\mathfrak{a}}
\newcommand{\fb}{\mathfrak{b}}
\newcommand{\fp}{\mathfrak{p}}
\newcommand{\fq}{\mathfrak{q}}
\newcommand{\fm}{\mathfrak{m}}
\newcommand{\fn}{\mathfrak{n}}
\newcommand{\FN}{\mathfrak{N}}
\newcommand{\FR}{\mathfrak{R}}
\newcommand{\set}[1]{\{#1\}}
\newcommand{\trv}{\set{1}}
\newcommand{\Aut}{\mathrm{Aut}}
\newcommand{\End}{\mathrm{End}}
\newcommand{\Ker}{\mathrm{Ker}}
\newcommand{\chr}{\mathrm{char}}
\newcommand{\Spec}{\mathrm{Spec}}
\newcommand{\Gr}{\mathrm{\bf Gr}}

\newcommand{\ord}{\operatorname{\mathrm{ord}}}
\newcommand{\ordp}{v}
\maketitle

\mbox{\hfill
  {\em For my adorable nephew, Sebastian Wayne Rusek, Born 6-16-2011.}
    \hfill\mbox{
      }}

\begin{abstract}

We study $A$-discriminants from a non-Archimedean point of view, refining earlier
work on the tropical discriminant. In particular, we study the case where $A$
is a collection of $n+m+1$ points in $\BZ^n$ in general position, and give an
algorithm to compute the image of the $A$-discriminant variety under the
non-Archimedean evaluation map. 
When $m=2$, our approach
yields tight lower and upper
bounds, of order quadratic in $n$. We also detail a Sage package for plotting
certain $p$-adic discriminant amoebae, and present explicit examples of point
sets yielding discriminant amoebae with extremal behavior.

\end{abstract}

\section{Introduction}
Amoebae --- the images of algebraic varieties under a valuation
map --- are of considerable interest in several complex variables,
tropical geometry, and arithmetic dynamics
\cite{gkz94,np09,br06}.
Furthermore, in addition to applications in mathematical physics
\cite{prt11}, amoebae have recently been used
to derive efficient algorithms in real algebraic geometry and
arithmetic geometry \cite{airr11, prt09}. In particular, the real part (and the
non-Archimedean rational part) of the complement of a discriminant amoeba
results in a new point of view in the classical study of discriminant
complements.

In this paper, we focus on the non-Archimedean amoebae of A-discriminants,
proving new complexity bounds on the topology of their closures.
In particular, we exhibit some unusual behavior differing from the complex
setting, and improve an earlier topological bound of Dickenstein, et.
al \cite{drrs07}.

Let $A$ be a generic collection of $n+m+1$ points in $\BZ^n$, with $m\ge0$. The
support of a polynomial is the collection of exponent vectors with nonzero
coefficients. The
$A$-discriminant of the family of polynomials over a field, $k$,  with support,
$A$, was introduced by Gelfand, Kapranov, and Zelevinsky in their book
\cite{gkz94}. They also discussed a simple parametric map called the Horn
uniformization of this $A$-discriminant. We will not discuss the original
parametrization further, but we will look at a dehomogenized form of it that
produces
the so-called reduced $A$-discriminant. Write $A=\{a_1,\dots,a_{n+m+1}\}$. Then
let $\hat{A}$ be the matrix 
\[ 
\left[ 
  \begin{array}{ccccc} 
    1&1&\cdots&1&1\\ 
    a_1&a_2&\cdots&a_{n+m}&a_{n+m+1}
  \end{array} 
\right].  
\] 
That is, we treat $A$ as a matrix whose column vectors are the $a_i$ and then we
let $\hat{A}$ be the same matrix with an extra row of all ones.
Let $B=\{\beta_{i}\}_{i=1}^{n+m+1}\subset\in\BZ^{m}$ be the right integer null space of
$\hat{A}$. Then the parametric plot for the reduced $A$-discriminant is a map
from $k^{m-1}$ to $k^m$ and has $\ell^{th}$ coordinate 
\[
\prod_{i=1}^{m+n+1}(\beta_{i,1}\lambda_1+\cdots+\beta_{i,m-1}\lambda_{m-1}+\beta_{i,m})^{\beta_{i,\ell}}.
\]
In our setting we will have $k$ as a non-Archimedian valuation field, so we can 
take the coordinate-wise
valuation. This will give us a piecewise linear object called the non-Archimedean
$A$-discriminant amoeba.

On the other hand, we have the real semi-algebra, $(\BR,\otimes,\oplus)$ where
$\otimes$ is standard addition and $\oplus$ gives the minimum of two real
numbers. Then given a polynomial in this semi-algebra, its zero set, a tropical
varieties, is defined
as the points where the graph is not differentiable. This too is a piecewise
linear object. It turns out that these two families of objects, non-Archimedean
amoebae and real semi-algebraic varieties have the same
combinatorial type\cite{PT05}.


Little has been written about explicitly and efficiently representing these
reduced
non-Archimedean amoebae. Kapranov's non-Archimedean theorem gives a construction
for non-Archimedean amoebae\cite{kap00}, but it requires constructing the
discriminant polynomial. This can be quite inefficient. For example the family
described in Example \ref{ex:rs} has a reduced discriminant polynomial with coefficients
with thousands of digits\cite{drrs07}. On the other hand, our software computes
the amoeba in seconds. Similarly Dickenstein, Feichtner, 
and Sturmfels fully
described the discriminant amoeba in the case where $k=\BC\{\{t\}\}$ in \cite{dfs}
and Rinc\'on built upon that same setting \cite{rincon11}. 
When $k$ is another field, such
as $\BC_p$, very little has been written in the direction.
In this paper, we make an effort to begin to close these
gaps. We begin with an explicit representation of the reduced
non-Archimedean discriminant amoeba as a collection of parametric tropical maps. 
Which lead themselves nicely to theorems.

Next we shift our attention to the special case $m=2$, where the
$A$-discriminant amoeba is a $1$-dimensional object in $2$-space. The process of
proving the previously mentioned reduction will further lead us to reductions in
the $m=2$ case. In the real case, the best known bound on the number of
components of the complement of these amoebae is $O(n^6)$ \cite{drrs07}, but our 
reductions naturally lead us to an upper bound that is $O(n^2)$. 
Furthermore, the method leads us to ways to produce extremal
examples. That is, we exhibit a family of $n$-variate $n+3$-nomials with
$O(n^2)$ connected components in the complement of the closure of their tropical
$A$-discriminant amoebae.
Finally, although it can be shown that the complex reduced $A$-discriminant is always
solid\cite{pst05} using these reductions we have found simple
polynomials (discussed in Example \ref{ex:rs}) that are not simply connected. A simple example is
  \[
  f(x,y,t)=tx^6+\frac{44}{31}ty^6-yt+y^6+\frac{44}{31}tx^3-x.
  \]
The associated $3$-adic $A$-discriminant is
\[
  \begin{picture}(130,125)
    \put(0,0){\includegraphics*[scale=0.5]{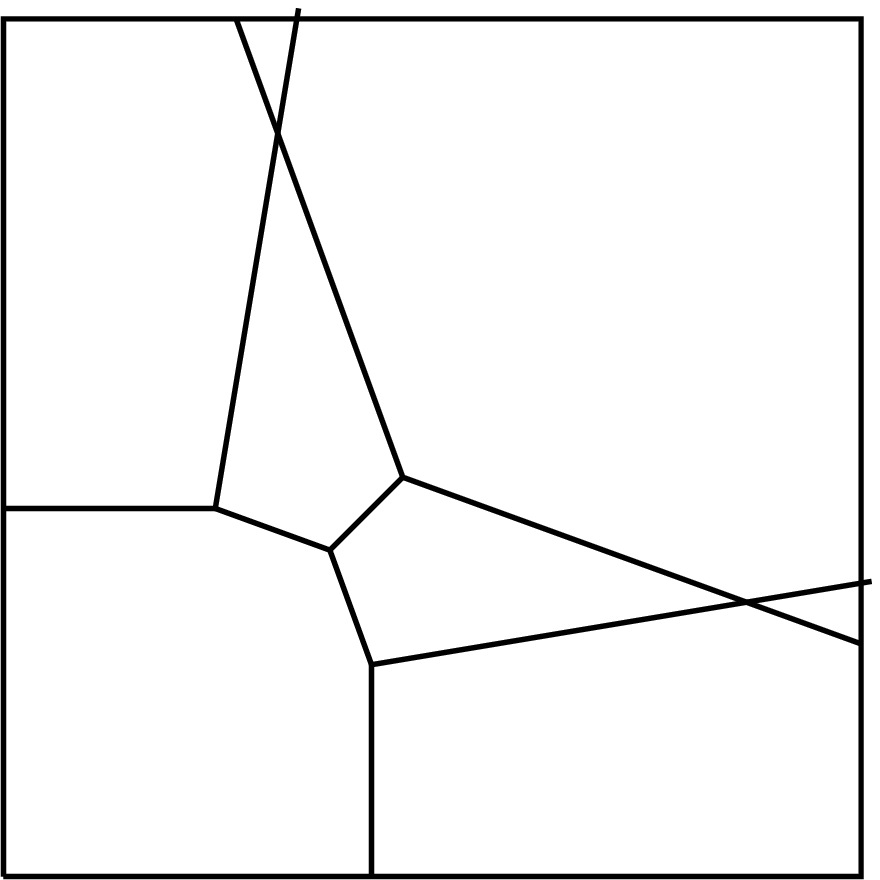}}
  \end{picture}
\]
which has two bounded connected components in its complement--two holes. This is
different than both the complex case and the standard tropical case
where $k=\BC\{ \{t\}\}$.

\section{Tropicalization of Parametric Non-Archimedean Maps}


We can construct
the reduced discriminant amoeba for the family of polynomials with support $A$,
which will be a $m-1$ dimensional surface in $m$-space.  We defined $A$ as a
collection of points, but we will abuse notation and also use it to represent
the $n$ by $(n+m+1)$ matrix whose columns are the points of $A$.
Then $\hat{A}\in\BZ^{(n+1)\times (n+m+1)}$ is the matrix $A$
with an extra row of all ones added. As in the introduction, we let
$B=\{\beta_{i,j}\in\BZ^{(n+m+1)\times m}$ be
the matrix whose columns form a basis of the right null-space of $\hat{A}$. Then
the map, $\phi:k^{m-1}\rightarrow\BR^{m}$ whose
$\ell^{th}$ coordinate is given by
\begin{equation}
  \sum_{i=1}^{n+m+1}\beta_{i,\ell}
\ordp(\beta_{i,1}x_1+\cdots+\beta_{i,m-1}x_{m-1}+\beta_{i,m}),
\end{equation}
where
$(\beta_i)$ are the column vectors of $B$, is a parametric form of the reduced
discriminant amoeba\cite{gkz94}. Throughout we will denote this function by $F$, the linear
forms $\beta_{i,1}x_1+\cdots+\beta_{i,m-1}x_{m-2}+\beta_{i,m}$ by $f_i$, and the
valuation of these linear maps by $F_i$. Often we will perform a linear
change of variables. This will often be denoted by a $*$ ($F^*$) or a
subscript ($F_I$) when necessary.

We want to better understand the structure of the image of this particular map,
$\phi$. 
Notice that the number in front of the valuation, $\beta_{i,\ell}$,also occurs
as a coefficient in the linear forms. We want to simplify $phi$ in such a way
that makes it easier to illustrate. To this end we wil look at a slightly more
general map and define a simplification. Thus the next two paragraphs will
deal with a slightly more general map instead.
Suppose instead that we have a map 
$G:k^{m-1}\rightarrow\BR^m$, a $p$-adic parametric plot whose $\ell^{th}$
coordinate is of the form
\begin{equation}\label{eqn:ord} \sum_{i=1}^{n+m+1}\gamma_{i,\ell}
  \ordp(\beta_{i,1}x_1+\cdots+\beta_{i,m-1}x_{m-1}+\beta_{i,m}), \end{equation}
where $\beta$ and $\gamma$ are $(n+m+1)\times m$ matrices with elements in
$k$.  Define the tropicalization of $G$ to be the parametric tropical map,
$\varphi_G:\BR^{m-1}\rightarrow\BR^m$, whose $\ell^{th}$ coordinate map is
\[\otimes_{i=1}^{n+m+1}(\ordp(\beta_{i,1})\otimes
r_1\oplus\cdots\oplus\ordp(\beta_{i,m-1})\otimes
r_{m-1}\oplus\ordp(\beta_{i,m}))^{\gamma_{i,\ell}},\] 
where the $\gamma_{i,j}$ are multiplied by the forms and we use the minimum in our
semi-algebra. 

Now if we let $g_i=\beta_{i,1}x_1+\cdots+\beta_{i,m-1}x_{m-1}+\beta_{i,m}$ then
we can rewrite the $\ell^{th}$ coordinate of $G$ as
\[
r\mapsto\sum_{i=1}^{n+m+1}\gamma_{i,\ell}
\ordp(g_i).
\]
This simplification will be a key piece in our paper. For now we will just state
how this is used in the process of simplifying the map $\phi$, but in the next
section we will fully prove the statements overviewed here.
Pick $I=\{i_1,\dots,i_{m-1}\}$ with $1\le i_1<\cdots<i_{m-1}\le n+m+1$ and such
that the zero set of $g_{i_1},\dots,g_{i_{m-1}}$ consists of a unique point in
$k$. Now, via Gauss-Jordan elimination
let $F_I$ be $F$ under the change of variables such that $f^*_{i_j}=x_j$.
Then we can write $F_I$ in a form similar to (\ref{eqn:ord}).
As we did with the original map, $F$, we can tropicalize our $F_I$ and get
$\phi_{F_I}$. Let $T$ be the collection of all such $I$. In the next
subsection,
we will prove
\begin{theorem}\label{thm:main}
Let $F$ be the parametric form of the non-Archimedean discriminant amoeba. Let $T$ be
the collection of all $I=\{i_1,\dots,i_{m-1}\}\subset\{1,\dots,n+m\}$ with 
$i_j<i_{j+1}$ for $j<m-1$. Let $F_I$ be $F$ under the linear change of variables
making $f_i=\lambda_i$ and $\phi_{F_I}$ the tropicalization of $F_I$. Then we
have
  \[
  \overline{F(k^{m-1})}=\bigcup_{I\in T}\phi_{F_I}(\BR^{m-1}).
\]
\end{theorem}

\subsection{Reducing non-Archimedean Amoebae to Parametric Tropical Functions}

A key difference between general non-Archimedean arithmetic and arithmetic in
the real semi-algebra is the ultrametric 
inequality. That is, when $a,b\in\BC_p$ then
$\ordp(a+b)\ge\min\{\ordp(a),\ordp(b)\}$,
whereas when $a$ and $b$ are elements of the real semi-algebra then $a\oplus
b=\min\{a,b\}$. For example, in $\BQ_p$, $\ordp(p+p^2)=\ordp(p)=1$, but
$\ordp(p+(p-1)p)=\ordp(p^2)=2$. In the latter sum the valuation of the sum is larger
than the valuation of the summands. We will call this {\it carrying}.
This carrying makes it appear that this non-Archimedean map and the
tropical maps may fail to have the same image because the carrying would
possibly cause
discontinuities in the image, whereas the image of the tropical map is
necessarily continuous. This section will show that this does not cause a
problem. That is, this section will prove theorem \ref{thm:main}.

Since the first proof in this subsection is long and notation heavy, we will begin
with an outline of the proof, so that the reader has an idea of the flow and
direction of the proof beforehand.
The proof begins by choosing an arbitrary parameter
$\lambda\in k^{m-1}$. In the next step we apply a change of variables on the
$f_i$. We will describe how the change of variables works.
The coefficients of the $f_i$ are taken from the rows of the $B$ matrix. There
are more rows than columns of the $B$ matrix. Now we may select a collection
$I=(i_1,\dots,i_{m-1})$ such that those rows of the matrix are linearly
independent. We then perform Gauss-Jordan elimination on the {\it columns} of
$B$ making each $i_j$ into a standard basis row vector. From this new matrix we
can construct a new collection of functions, $\{f^*_1,\dots,f^*_{n+m+1}\}$ with
$f^*_{i_j}=x_j$. These functions represent the original collection after a
linear change of variables. This is the linear change of variables we will use,
and the linear change of variables that results from the proof. When we apply
this change of variables to $F$, we will call the resulting map $F_I$. When the
correct $I$ (based on $\lambda$) is chosen the map
will have the
property that it
exhibits a simple parameter $z_i=p^{\ell_i}$ that approximates $F(\lambda)$.
Moreover, we show that we can pick the approximating $z_i$ in a special 
way such that no non-Archimedean carrying happens, on a dense subset of the
domain. Since there is little carrying 
then we can
replace the map $F^*$ with its polarization and the parameter $p^{\ell_i}$
with $\ell_i$. Hence we can approximate $F(\lambda)$ using one of the
$\phi_{F_I}$ as defined earlier. Since the $\phi_{F_I}$ are closed maps
then we will have one containment. The other containment will follow quite easily.

To further illustrate this change of variables we will go through an example.
The example will explicitly construct a couple of the $F_I$ and $\phi_{F_I}$ and will
illustrate the difference between the images of these maps.

\begin{example}
  This example will take place in 2-space. The same idea presents itself in
  higher dimensions, but it is easier to grasp here. We will work in $\BQ_2$.
  In this simple example $n=0$ and $m=2$.
  Let
  \[
    f_1=\lambda-1,~f_2=\lambda-13,~f_3=\lambda-25.
  \]
  Hence we have
  \[
  F(\lambda)=(\ordp(\lambda-1)+\ordp(\lambda-13)+\ordp(\lambda-25),
  -\ordp(\lambda-1)-13\ordp(\lambda-13)-25\ordp(\lambda-25)).
  \]
  Now $F_{\{1\}}$ requires $f^*_1=\lambda$, so we would have
  \[
  f^*_1=\lambda,~f^*_2=\lambda-12,~f^*_3=\lambda-24.
  \]
  This would give us
  \[
  \phi_{F_{\{1\}}}(r)=(
  (r)\otimes(r\oplus2)\otimes(r\oplus3),
  -1(r)\otimes-13(r\oplus2)\otimes-25(r\oplus3)).
  \]
  (Note that the linear forms change but the coefficients in front of them do
  not.)
  Now when $\ell\ne2,3$ then $F^*_1(2^\ell)=\ell$,
  $F^*_2(2^\ell)=\min\{\ell,2\}$, $F^*_3(2^\ell)=\min\{\ell,3\}$, so when
  $\ell\ne2,3$ we have
  \[
  F^*(2^\ell)=(\ell+\min\{\ell,2\}+\min\{\ell,3\},
  -1\ell-13\min\{\ell,2\}-25\min\{\ell,3\}).
  \]
  But $F^*$ 
  has a discontinuity at $\ell=3$. 
  That is,
  $F^*_3(2^3)=\ordp(16)=4$, but is less than or equal to $3$ everywhere else, so
  the point $F^*(2^3)=(9,-129)$ is an isolated point when restricting to $2^\ell$.
  (A similar problem results from $\ell=2$.) On
  the other hand $\phi_{F_{\{1\}}}$ does not have this discontinuity.
  The map $\phi_{F_{\{1\}}}$ does not have
  a carry and $\phi_{F_{\{1\}}}(3)=(8,-78)$, so this value, $(9,-129)$, appears to be lost.  
  Furthermore, it appears this point, $(9,-129)$, will be an isolated point in the
  image and is not part of the images of the $\phi_{F_{\{i\}}}$.
  What has happened here? 
  When we instead apply the change of variables used in making
  $\phi_{F_{\{3\}}}$ 
   we see what happens. That is
  \[f'_1=\lambda+24,~f'_2=\lambda+12,~f'_3=\lambda\] 
  and 
  \[
  \phi_{F_{\{2\}}}(r)=(
  (r\oplus3)\otimes(r\oplus2)\otimes(r),
  -1(r\oplus3)-13\otimes(r\oplus2)\otimes-25(r)).
  \]
  Now when $r=4$ we have $\phi_{F_{\{2\}}}(4)=(9,-129)$. Thus we see an
  illustration of how though the $F_{\{i\}}$ do not have the exact same image as
  their corresponding $\phi_{F_{\{i\}}}$, but the collection still contains the
  desired values.

\end{example}

We will now prove that $F(\lambda)$ can be approximated as desired.
\begin{theorem}\label{thm:red1}
  Let $k$ be an algebraically closed complete valuation field with
  $\BQ\subseteq\ord(k)\subseteq\BR$.
  Pick any $\lambda=(\lambda_1,\dots,\lambda_{m-1})\in k^{m-1}$ where $F(\lambda)$ is
  well-defined. That is, $f_i(\lambda)\ne0$ for all $i$. Let $\omega\in k$ be
  any element with $\ord(k)=1$ 
  (we would naturally select $\omega=t$ and $\omega=p$, for
  $k=\BC\{\{t\}\}$ and $k=\BC_p$, respectively). 
  For any
  $\varepsilon>0$
  there is $I=(i_1,\dots,i_{m-1})\in\{1,\dots,n+m+1\}^{m-1}$  and an $F^*$
  that is $F$ under a linear change of variables (Gauss-Jordan Elimination) such that
  $F^*_{i_j}(x)=x_j$ such that there are $\ell_1,\dots,\ell_{m-1}\in\BQ\setminus\BZ$
  with $\ell_i-\ell_j\not\in\BZ$ for $i\ne j$ with
  $|F^*(\omega^{\ell_1},\dots,\omega^{\ell_{m-1}})-F(\lambda)|<\varepsilon$.
  \begin{proof}
    Pick $\lambda\in k^{m-1}$ as described in the hypothesis and choose any
    $0<\varepsilon<1$.  For any $z\in k$ let
    $\lambda_z=(\lambda_1+z,\lambda_2,\dots,\lambda_{m-1})$. Hence
    $f_i(\lambda_z)=f_i(\lambda)+\beta_{k,1}z$.  Now I claim that for each $i$ with
    $\beta_{i,1}\ne0$, there is an $N_i$ such that for any $z\in k$ with
    $\ordp(z)>N_i$ then $F_i(\lambda_z)=F_i(\lambda)$. On the other hand, for
    this same $i$, if $\ordp(z)<N_i$ then
    $F_i(\lambda_z)=\ordp(z)+\ordp(\beta_{i,1})$.  Indeed, if for any $i$, we pick
    $N_i=F_i(\lambda)-\ordp(\beta_{i,1})$ we
    have
    \begin{align*}
      F_i(\lambda_z)&\ge\min\{\ordp(f_i(\lambda)),\ordp(z)+\ordp(\beta_{i,1})\}\\
      &\ge\min\{N_i+\ordp(\beta_{i,1}),\ordp(z)+\ordp(\beta_{i,1})\}.  \end{align*}
    Now when $\ordp(z)\ne N_i$ then this ultrametric inequality becomes an
    equality. Hence when $\ordp(z)<N_i$ then we have
    $F_i(\lambda_z)=\ordp(z)+\ordp(\beta_{i,1})$ and when $\ordp(z)>N_i$ we have
    $F_i(\lambda_z)=N_i+v(\beta_{i,1})=F_i(\lambda)$. Now at least one $i$ should exist with $\beta_{i,1}\ne0$ or
    $x_1$ plays no role in any of our equations. Without loss of generality
    assume that $1$ be the (an) index associated with the maximum
    such $N_i$ and let $N_{1}$ be the
    relevant value. This will be our $i_1$.

    To show that this works as our $i_1$,
    we pick $\ell\in\BQ$ such that $0<N_{1}-\ell<\varepsilon$.  Now if for
    all $i$, we
    let $f'_i$ be $f_i$ under the change of variables $x_1\mapsto
    x_1-\frac{\beta_{1,2}}{\beta_{1,1}}x_2-\cdots-\frac{\beta_{1,m-1}}{\beta_{1,1}}x_{m-1}
    -\frac{\beta_{1,m}}{\beta_{1,1}}$. Then we have
    \begin{align*} 
      f'_i(\omega^\ell,\lambda_2,\dots,\lambda_{m-1})&=
      \left(\beta_{i,1}\omega^\ell-\frac{\beta_{i,1}\beta_{1,2}}{\beta_{1,1}}\lambda_2-\cdots-\frac{\beta_{i,1}\beta_{1,m-1}}{\beta_{1,1}}\lambda_{m-1}
      -\frac{\beta_{i,1}\beta_{1,m}}{\beta_{1,1}}\right)\\
      &\phantom{=}+\beta_{i,2}\lambda_2+\cdots+\beta_{i,m-1}\lambda_{m-1}+\beta_{i,m}\\ &=
      \beta_{i,1}\omega^\ell-\left(\frac{\beta_{i,1}\beta_{1,1}}{\beta_{1,1}}\lambda_1
      +\frac{\beta_{i,1}\beta_{1,2}}{\beta_{1,1}}\lambda_2+\cdots+\frac{\beta_{i,1}\beta_{1,m-1}}{\beta_{1,1}}\lambda_{m-1}
      +\frac{\beta_{i,1}\beta_{1,m}}{\beta_{1,1}}\right)\\ &\quad+\left(\beta_{i,1}\lambda_1+
      \beta_{i,2}\lambda_2+\cdots+\beta_{i,m-1}\lambda_{m-1}+\beta_{i,m}\right)\\
      &=\beta_{i,1}\omega^\ell-\frac{\beta_{i,1}}{\beta_{1,1}}f_1+f_i 
    \end{align*}
    We desire to know the valuation of
    $f'_i(\omega^\ell,\lambda_2,\dots,\lambda_{m-1})$. To this end, we look at the
    valuation of the three pieces in this sum.
    The first pieces gives $\ordp(\beta_{i,1}\omega^\ell)=\ordp(\beta_{i,1})+\ell$. The
    second,
    $\ordp(\frac{\beta_{i,1}}{\beta_{1,1}}f_1) =
    N_1+\ordp(\beta_{1,1})+\ordp(\beta_{i,1})-\ordp(\beta_{1,1}) = N_1+\ordp(\beta_{i,1})$,
    since we have chosen $N_1$ such that $\ordp(f_1)=N_1+\ordp(\beta_{1,1})$. Similarly, the third piece
    gives
    $\ordp(f_i)=N_i+\ordp(\beta_{i,1})$. By the ultrametric inequality, if these
    three items have different valuation then their sum has valuation that is
    the minimum of the three.
    Therefore if $N_i<N_1$ then $N_i<\ell$ and
    $F'_i=F_i$. Otherwise, if $N_i= N_1$ (it cannot be larger) then $\ell<N_1= N_i$, which tells us
    that $F'_i=\ell+\ordp(\beta_{i,1})$ and so $|F'_i-F_i|\le\varepsilon$.
    We see that the latter two pieces, $f_i$ and
    $\frac{\beta_{i,1}}{\beta_{1,1}}f_1$ have the same coefficient on $x_1$ so the
    difference is independent of $x_1$. That is, the other $F'_i$ are independent
    of $x_1$, and changing
    the value of $\ell$ doesn't affect their difference. Hence we can make
    another change of variables sending $x_1$ to $\frac{x_1}{\beta_{1,1}}$ and use
    an $\ell$ that is $\ordp{\beta_1}$ larger. This makes $f'_1=x_1$, as
    desired. Then $i_1=1$ and $\ell_1=\ell$. 
    We can go through this procedure again with $x_2$ and the modified
    collection $\{f'_i\}$. It is clear that the newly selected $i_2$
    will not be $1$, because we have made $f'_{1}$ depend only
    on $x_1$. We can continue this iteratively through all the variables. Each
    time one of the $f_i$ will necessarily be chosen or the parametric function
    is under determined.
    Selecting $i_1,\dots,i_{m-1}$ and
    $\ell_1,\dots,\ell_{m-1}$
    and a final collection of $f^*_{i_{i_1}},\dots,f^*_{i_{m-1}}$. Then the final linear
    forms $\{f^*_i\}$ clearly do what we want and
    $F^*(p^{\ell_1},\dots,p^{\ell_{m-1}})$ approximates $F(\lambda)$ to within
    $b\varepsilon$ where $b$ depends only on the original matrix, $B$.
  \end{proof} 
\end{theorem}

It is again worth noting that the linear change of variables used was Gauss-Jordan
elimination on the $f_{i_1},\dots,f_{i_{m-1}}$. This means for a given
collection, if we sort the $i_1,\dots,i_{m-1}$ there is a well-defined change of
variables to use together with a finite collection of choices on the whole.
Now theorem \ref{thm:main} will follow quite easily:

\begin{theeorem}
  Let $T=\{I=(i_1,i_2,\dots,i_{m-1})~~|~~0\le i_1<i_2<\cdots<i_{m-1}\le
  n+m+1~with~\#Z(\{f_{i_j}\})<\infty\}$. Then
  \[\overline{F(k^{m-1})}=\bigcup_{I\in T}\varphi_{F_I}(\BR^{m-1}).\]
  \begin{proof}
    Theorem \ref{thm:red1} says that for any $z\in k^{m-1}$ there are an $I$ as
    described above and an
    $\ell\in\BQ^{m-1}$ such that $|\varphi_{F_I}(\ell)-F(z)|<\varepsilon$. Thus
    since the $\varphi_{F_I}$ are closed maps and there are finitely many of
    them then $\overline{F(k^{m-1})}\subset\bigcup\varphi_{F_I}(\BR^{m-1})$.

    Now it suffices to show that $F(k^{m-1})$ is dense in
    $\varphi_I(\BR^{m-1})$. Pick any $I$ with $Z(f_I)$ finite
    and any $\ell\in\BQ^{m-1}$. We may assume that
    $\ell_i-\ell_j\notin\BZ$ for $i\ne j$ and $\ell_i\notin\BZ$. Clearly the
    collection of all such $\ell_1,\dots,\ell_{m-1}$ is still dense in
    $\BR^{m-1}$. Let $F^*$ be $F$
    under the linear change of variables making $f_{i_j}=x_j$, as described in
    the previous theorem.
    Then it is clear that 
    \[F^*_i(\omega^{\ell_1},\dots,\omega^{\ell_{m-1}})=\min_j\{\ordp(a^*_{i,j}\omega^{\ell_j})\},\]
    because no carrying can occur since $\ell_i-\ell_j\notin\BZ$.
    Therefore
    $F^*(\omega^\ell_1,\dots,\omega^{\ell_{m-1}})=\varphi_I(\ell_1,\dots,\ell_{m-1})$ and so
    $\overline{F(k^{m-1})}\supset\bigcup\varphi_{F_I}(\BR^{m-1})$ as desired.
  \end{proof}
\end{theeorem}

\section{The Case $m=2$ and $k=\BQ_p$}\label{sec:m2}

When we apply the change of variables in the main theorem we are essentially
approximating the zero of a collection of $m-1$ linear forms. We see this
because as $\ell_i$ goes to infinity the $p$-adic parameter $p^{\ell_i}$ goes to
zero and hence so do the associated linear forms. In the special case where
$m=2$ this is exactly what we are doing because there is only one parameter. 
To simplify notation in this section let $a$ and $b$
in $\BR^{n+m-1}$ be the columns of $B$, and let $\phi_i$ be the polarization of
$F$ after applying the variable change making $f_i=\lambda$.  Given a particular
$i$ where $a_i\ne0$, we have $f_i=a_i\lambda+b_i$ and the desired change of
variables producing $\phi_i$ is
$\lambda\mapsto \frac{p^{\ell}-b_i}{a_i}$. As $\ell$ goes to infinity, $\lambda$
approaches $\frac{-b_i}{a_i}$, namely, the zero of $f_i$.

Let $z_i=\frac{-b_i}{a_i}$. Now if $z_i\equiv z_j\mod p^q$ then
$\phi_i(\ell)=\phi_j(\ell)$ for $\ell\le q$ because we are using strict minimum.
Remember that $F_i:=\ordp(f_i)$. When $z_i\equiv z_j\mod p^q$ then
$F_i(z_j+p^\ell)=\ell+\ordp(a_i)$ for $\ell<q$, because $F_i(z_j)\ge
q+\ordp(a_i)>\ell+\ordp(a_i)$. On the other hand for $\ell>q$ we have that
$F_i(z_j+p^\ell)$ is the constant $q+\ordp(a_i)$. These facts mean that for a
particular $i$, $\phi_i$ is linear (or constant) everywhere except at the
collection of points $\ordp(z_i-z_j)$ for $j\ne i$. In the language of tropical
geometry, we are saying that $V=\{\ordp(z_i-z_j)\}_j$ is the tropical hypersurface
of the parametric plot $\phi_i$.  

This tells us that there is overlap between the various $\phi_i$ for various
values of the parameter $\ell$. With this in mind we will create a tree mapping
out the possible differences. A nonzero element $a\in\BQ_p^\times$ with $\ordp(a)=j$
can be written $a=\sum_{i\ge j}a_ip^i$ with $a_i\in\BZ/p\BZ$ and $a_j\ne0$. For
a given $a=\sum_{i\ge j}a_ip^i$, we will call $a_i$ the {\it digit} of $a$ at
$i$. If $i<j$ then the digit of $a$ at $i$ is defined to be zero. Thus we can
construct a tree expressing the relationships between these elements. The head
node represents the smallest value in $\{\ordp(z_i-z_q)~|~i \ne q\}$. If $j$ is
this value then the tree will have a branch for each distinct digit at $j$ among
all the $z_i$. Now each branch represents a different digit at $j$, and we
associate the elements with that digit at $j$ to the branch associated to that
digit at $j$. For a given branch we follow the same procedure, but only with
the $z_i$ associated to that branch. We repeat this iteratively and eventually a
branch will have only a single element associated to it. We add one more node to
the end of this branch and label it with the given $z_i$.

Now these two paragraphs together tell us that a given branch, between two
non-leaf nodes, represents a line segment in the plot of the amoeba. We also 
see that every element
associated to that branch contributes the value of the coefficient multiplied by
its linear form to the slope of the branch. This is because these are precisely
the forms producing the parameter plus a constant and the other elements are the
ones producing a simple constant. The branches connected to leaf nodes, on the 
other hand, still have slope that is the value of the number multiplied by the
associated linear form, but they represent a ray, because now the parameter,
$\ell$, goes to infinity. Finally, there is possibly one more ray. When the
parameter, $\ell$, goes to negative infinity, then we can have another ray. This
only happens when the coefficients on the linear forms don't add to zero. Though
we are assuming the $a_i$ and the $b_i$ add to zero, we can get a nonzero sum in
the case that one of the $a_i$ is zero, because then the relevant $b_i$ will not
be used in the sum. We give an example to illustrate what we've mentioned.

\begin{example}\label{ex:rs}
  Consider the support 
  \[
  \left[
  \begin{array}{cccccc}
    6&0&0&0&3&1\\
    0&3&1&6&0&0\\
    1&1&1&0&0&0
  \end{array}
  \right].
  \]
  This is the support of the so-called Rusek-Shih example 
  \[
  f(x,y,t)=tx^6+\frac{44}{31}ty^6-yt+y^6+\frac{44}{31}tx^3-x
  \]
  from \cite{drrs07}.
  We then choose our $B$ matrix to be the transpose of the following:
  \[
  \left[
  \begin{array}{rrrrrr}
    -2&35&-33&-12&0&12\\
    2&-11&9&4&-4&0
  \end{array}
  \right]
  \]
Here is an example tree. If our $z_i$ are $1,11/35,3/11,1/3,$ and $0$ then we
have the following $3$-adic expansions starting from index $-1$.
\[
\begin{array}{lcr}
  \alpha_0 =&\frac{1}{3}   &= [1,0,0,0,\dots]\\
  \alpha_1 =&\frac{11}{35} &= [0,1,2,2,\dots]\\
  \alpha_2 =&1             &= [0,1,0,0,\dots]\\
  \alpha_3 =&\frac{3}{11}  &= [0,0,2,1,\dots]\\
  \alpha_4 =&0             &= [0,0,0,0,\dots]
\end{array}
\]
These elements would then be put into the previously described tree as
  \begin{center}
  \begin{picture}(240, 85)
    \put(0,5){\includegraphics{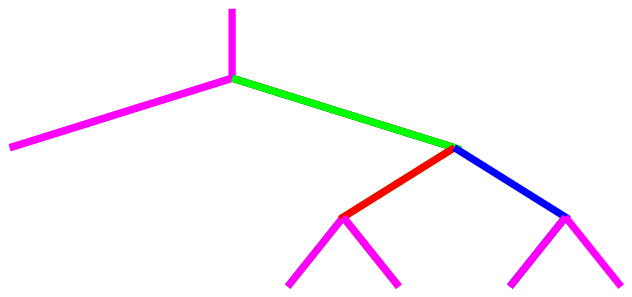}}
    \put(53,41){$\alpha_0$}
    \put(133,1){$\alpha_1$}
    \put(163,1){$\alpha_2$}
    \put(197,1){$\alpha_3$}
    \put(227,1){$\alpha_4$}
  \end{picture}
\end{center}
Notice that the elements to the left of the first node on the right begin with
$[0,1]$ while the ones to the right begin with $[0,0]$ and similar relations can
be seen below the other nodes.
Each branch of this tree represents a segment (or ray, for the leaves) of
constant slope. In particular, with this example, when we plot the amoeba we get
\[
  \begin{picture}(130,125)
    \put(0,0){\includegraphics*[scale=0.5]{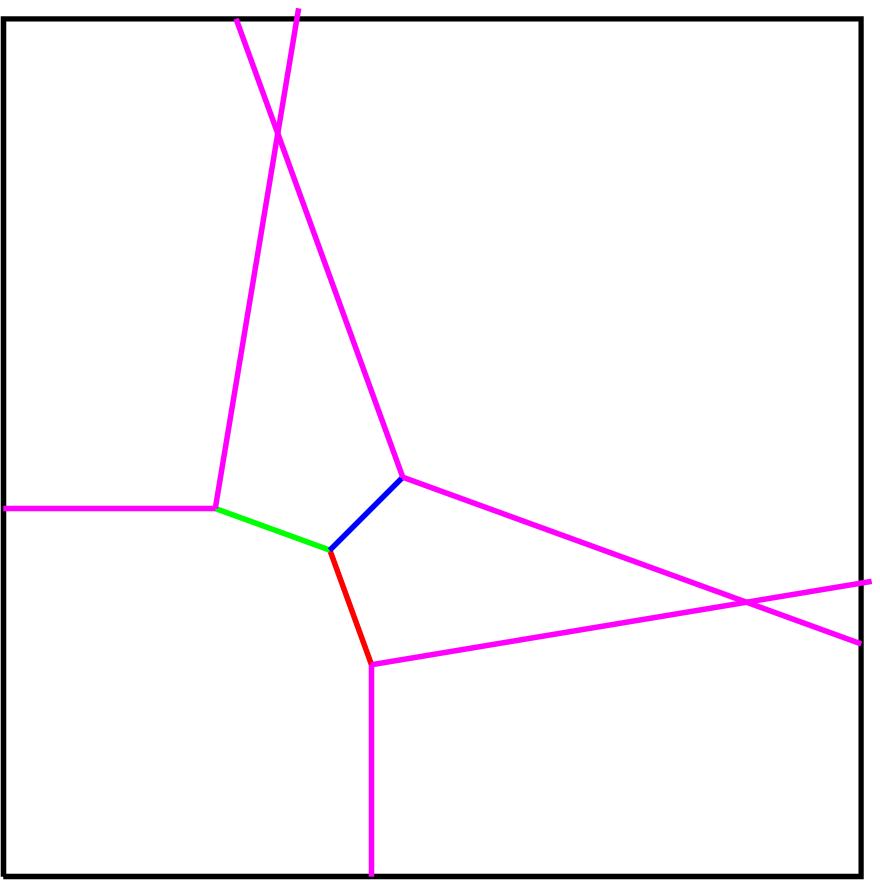}}
  \end{picture}
\]
The colors are there to help indicate which branch corresponds to which segment
or ray. Notice the node in the middle of the graph that is not connected to any
ray. This corresponds to the node earlier mentioned in the graph, because it has
no leafs as direct decendants. Also, there is an extra ray. This accounts for 
$\ell$ approaching
negative infinity and the fact that one of the last two linear forms has no zero.
\end{example}

This $p$-nary tree we have constructed has $n+3$ leaves, when none of the $a_i$
are zero, otherwise we include the extra ``leaf'' on top like our previous
example. It is a basic fact of graph theory that such a tree will have
no more than $2n+4$ branches, when $n>1$. As explained earlier each branch
represents a ray or a segment in the graph of the $p$-adic amoeba. That is, each
branch represents a straight piece of the amoeba. This fact will give us an
upper bound on the number of connected components of the complement of the
amoeba. If we replace each segment or ray with an entire line then we have a
line arrangement. It is then well known that such a line arrangement has no more
than $\binom{2n+4}{2}+\binom{2n+4}{1}+\binom{2n+4}{0}=2n^2+9n+11$ chambers in
its complement. Therefore the $p$-adic amoeba also has no more than
$2n^2+9n+11$ complement components when none of the $a_i$ are zero. On the other
hand, when an $a_i$ is zero we have fewer nodes and hence $2$ fewer branches,
but we have one extra ray accounting for $\ell\rightarrow-\infty$. Thus the same
bound still applies.

It seems unlikely that this is a strict upper bound, but it is not too difficult
to generate examples that have $O(n^2)$ examples. That is, we can show that the
bound is asymptotically tight. We write this as a theorem and in 
section \ref{sec:extreme} we will look at a family of examples exhibiting the asymptotic bound.
\begin{theorem} 
  The closure of reduced $A$-discriminant of $n+3$ points in general position
  has no more than $2n^2+9n+11$ complement components. 
  Moreover, Section \ref{sec:extreme} evinces supports $A_n$ in
  $\BZ^n$ with cardinality $n+3$ and primes $p_n$ such that the $p_n$-adic
  discriminant amoeba has quadratically (quadratic in $n$) many connected components in
  its complement.
\end{theorem}

%
%
%

\section{Extremal $p$-adic Family}\label{sec:extreme}

  We will construct a family of $A$ matrices admitting quadratically many
  complement components. We begin by constructing a $B$ matrix that satisfies
  our requirements and then
  work backwards from there to get the $A$ matrix.
  Let $p$ be a prime number and let $k$ be any integer larger than 2. We define
  $D\in\BZ^{2\times (2k+2)}$, by $D_{2,2i}=p^{i-1}$ and
  $D_{2,2i+1}=-p^{i-1}$ for $i=1,\dots,k+1$ and $D_{1,1}=D_{1,2}=-k$ and
  $D_{1,j}=1$ for $j>2$. That is, $D$ has the form:
  \[\left[
  \begin{array}{ccccccccc}
    -k&-k&1&1&1&1&\cdots&1&1\\
    -1&1&-p&p&-p^2&p^2&\cdots&-p^{k}&p^{k}\\
  \end{array}
  \right].\]
  Our $B$ matrix would be the transpose of $D$.
  Now the zeros of our linear forms are $\pm\frac{1}{k},\pm p,\pm p^2,\dots,\pm
  p^n$. The $p$-adic order of the first two elements is less or equal to $0$ and
  the others are their respective exponents on $p$. Therefore the associated tree
  is rather easy to form. For example, for $k=3$ and $p=2$ we have:
  \begin{center}
    \begin{picture}(240, 100)
      \put(0,5){\includegraphics{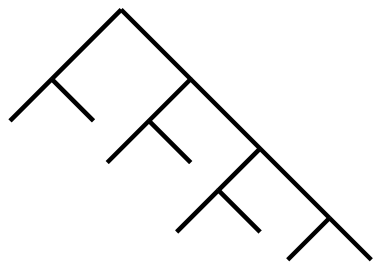}}
      \put(132,32){$-p$}
      \put(115,32){$p$}
      \put(152,12){$-p^2$}
      \put(135,12){$p^2$}
      \put(182,2){$-p^3$}
      \put(165,2){$p^3$}
      \put(99,40){$-\frac{1}{3}$}
      \put(82,40){$\frac{1}{3}$}
    \end{picture}
  \end{center}
  For larger $k$ the tree extends further to the right, whereas for larger $p$
  the leaves for $p^i$ and $-p^i$ branch directly from the main branch on the
  right rather than having their own mutual branch first. This is because the $i^{th}$
  digits for $p^i$ and $-p^i$ is the same only when $p=2$. 
  For example,
  $9=0\cdot3+1\cdot3^2,3=1\cdot3,$ and $-3=2\cdot3+3\cdot3^2+3\cdot3^3+\cdots$,
  while $4=0\cdot2+1\cdot2^2,2=1\cdot2,$ and $-2=1\cdot2+1\cdot2^2+\cdots$.
  That is, $\ord_3(3-(-3))=1=\ord_3(3-0)$, while $\ord_2(2-(-2))=2$.

  The slope of a branch is the sum of the $(a_i,b_i)$ of the $z_i$ associated to
  that branch.
  Therefore any non-leaf branch on the right has a slope of the form $(m,0)$, because each
  $p^i$ will cancel out with its matching $-p^i$, but the $1$'s in the first
  coordinate will add. Now at the branch point between $p^i$ and $-p^i$ we have
  rays in the direction $(1,p^i)$ and $(1,-p^i)$. That is, we have a line in the
  positive $x$-direction with rays emanating with slopes $\pm
  p^{-i}$. Furthermore each successive
  ray in the direction $(1,p^i)$ is further along the $x$-axis than the
  previous one because its associated branch in the tree splits further along
  the main branch. Therefore, 
  because the slope is
  less steep, the ray for $p^i$ (resp.~$-p^i$) intersections the ray for $p^j$
  (resp.~$-p^j$) for all $j>i$. The points from which these rays are emanating
  are independent of $p$. Thus for each $k$, a $p$ can be chosen assuring the
  intersections are non-degenerate. That is, for large enough $p$ the $p^i$ ray
  intersects the $p^j$ with a smaller $x$-coordinate than the starting
  position of the $p^{j+1}$ ray. Hence the rays above the $x$-axis give a
  line arrangement with
  at least $\binom{k}{2}+\binom{k}{1}+1$ components. Similarly the
  rays below the $x$-axis give the same number of components, except one of
  these components on the right is already accounted for in the previous count.
  This gives us at least $k^2+k+1$ components in the complement of the amoeba.
  As a visual example, here is the relevant part of the 
  discriminant amoeba for $p=3$ and $k=3$. You can also see the far right
  chamber that is not cut into two.
  \begin{center}
    \begin{picture}(200, 100)
      \put(0,0){\includegraphics*{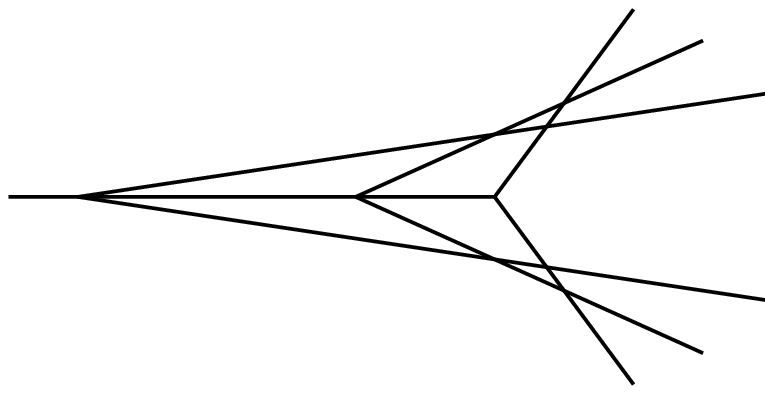}}
    \end{picture}
  \end{center}
  Now $k=\frac{n+m-2}{2}$. Hence the number of components is quadratic in $n$.

  Now constructing an $A$ matrix to accompany such a $B$ matrix is not hard.
  First we find the null space, $N$, of $D$. It will be a $2k$ by $2k+2$ matrix.
  For $i=1,\dots,k$, the odd rows will be $N(2i-1,1)=-kp^i+1$ and
  $N(2i-1,2)=kp^i+1$ and $N(2i-1,2i+1)=2k$, all other coordinates of that row
  being zero. Similarly, for $i=1,\dots,k$, the even rows will be
  $N(2i,1)=kp^i+1$, $N(2i,2)=-kp^i+1$, and $N(2i,2i+2)=2k$, again all other
  coordinates being zero.  It is clear that the rows of $N$ are linearly
  independent and it is the proper dimension. A small bit of arithmetic verifies
  that the rows of $N$ are orthogonal to the columns of $B$. Finally we can
  remove any row of $N$ to get our desired $A$ matrix. We remove one row because
  $B$ should be the null space of $\hat{A}$ rather than the null space of $A$.
  Therefore $A$ can have the form
  \[
  \left[
  \begin{array}{rrccccccccc}
    -kp  +1&kp  +1&2k&0&0&0&\cdots&0&0&0&0\\
    kp  +1&-kp  +1&0&2k&0&0&\cdots&0&0&0&0\\
    -kp^2+1&kp^2+1&0&0&2k&0&\cdots&0&0&0&0\\
    kp^2+1&-kp^2+1&0&0&0&2k&\cdots&0&0&0&0\\
        &     & & &  & &\cdots& & & &\\
        -kp^{k-1}+1&kp^{k-1}+1&0&0&0&0 &\cdots&2k&0&0&0\\
    kp^{k-1}+1&-kp^{k-1}+1&0&0&0&0 &\cdots&0&2k&0&0\\
    -kp^k+1&kp^k+1&0&0&0&0 &\cdots&0&0&2k&0\\
  \end{array}
  \right]
  \]

\section{Sage Code}

When $m=2$ this easily lends itself to a simple algorithm. 
A Sage package can be found at
``http://math.tamu.edu/\~{}krusek/pamoeba.sage''. One inputs the $B$ matrix
and the $p$-adic field to use and line objects are returned that represent the
segments and rays of the $p$-adic amoeba. The code itself constructs the tree
described in section \ref{sec:m2} then for each branch it creates a sage line
object representing the image of that branch.
Here is a short code snippet
to plot the extremal example from the previous section with $k=3$ and $p=3$.
\begin{verbatim}
load "pamoeba.sage"
B = MatrixSpace(ZZ,2,8)([[-3,-3,1,1,1,1,1,1],[-1,1,-3,3,-9,9,-27,27]])
B = B.transpose()
K = Qp(3)
lns = getamoeba2(B,K)
show(sum(lns))
\end{verbatim}

\bibliographystyle{plain}
\bibliography{plot-nonarch-amoeba}

\begin{thebibliography}{10}

\bibitem{airr11}
Martin Avenda\~no, Ashraf Ibrahim, J.~Maurice Rojas, and Korben Rusek.
\newblock Faster p-adic feasibility for certain multivariate sparse
  polynomials.
\newblock {\em Journal of Symbolic Computation}, 2011.
\newblock To appear.

\bibitem{br06}
Matthew Baker and Robert Rumely.
\newblock Equidistribution of small points, rational dynamics, and potential
  theory.
\newblock {\em Ann. Inst. Fourier (Grenoble)}, 56(3):625--688, 2006.

\bibitem{dfs}
Alicia Dickenstein, Eva~Maria Feichtner, and Bernd Sturmfels.
\newblock Tropical discriminants.
\newblock {\em The Journal of the American Mathematical Society},
  20:1111--1133, 2007.

\bibitem{drrs07}
Alicia Dickenstein, J.~Maurice Rojas, Korben Rusek, and Justin Shih.
\newblock Extremal real algebraic geometry and a-discriminants.
\newblock {\em Moscow Mathematical Journal}, 7(3):425--452, 2007.

\bibitem{rincon11}
Rinc\'on Filipe.
\newblock Computing tropical linear space.
\newblock {\em To appear in the Journal of Symbolic Computation}, page~15,
  2012.

\bibitem{gkz94}
Israel~Moseyevitch Gel'fand, Misha~M. Kapranov, and Andrei~V Zelevinsky.
\newblock {\em Discriminants, Resultants and Multidimensional Determinants}.
\newblock Mathematics: Theory \& Applications. Birkh\"auser, Boston, 1994.

\bibitem{kap00}
Mikhail~M. Kapranov.
\newblock Amoebas over non-archimedean fields.
\newblock {\em University of Toronto}, 2000.

\bibitem{np09}
Lisa Nilson and Mikael Passare.
\newblock Discriminant coamoebas in dimension two.
\newblock {\em pre-print}, page~17, 2009.

\bibitem{prt11}
Mikael Passare, Dmitry Pochekutov, and August Tsikh.
\newblock Amoebas of complex hypersurfaces in statistical thermodynamics.
\newblock {\em pre-print}, page~18, 2011.

\bibitem{pst05}
Mikael Passare, Timur Sadykov, and August Tsikh.
\newblock Singularities of hypergeometric functions in several variables.
\newblock {\em Compos. Math.}, 141(3):787--810, 2005.

\bibitem{PT05}
Mikael Passare and August Tsikh.
\newblock Amoebas: Their spines and their contours.
\newblock {\em Contemporary Math}, 377:275--288, 2005.

\bibitem{prt09}
Philippe Pebay, J.~Maurice Rojas, and David Thompson.
\newblock Optimizing n-variate (n+k)-nomials for small k.
\newblock {\em Theoretical Computer Science, Symbolic-Numeric Computation},
  412(16):1457--1469, 2011.
\newblock 2009 Special Issue.

\end{thebibliography}

\end{document}